\documentclass[12pt,reqno]{amsart} \usepackage{setspace}
\usepackage{amsmath, amssymb, amsthm}
\usepackage{amsmath}
\usepackage{amscd,amsthm,amsfonts,amsopn,amssymb,mathrsfs, latexsym}
\usepackage{epsfig,hyperref}
\numberwithin{equation}{section}
\theoremstyle{plain}
\newtheorem{theorem}{Theorem}

\theoremstyle{remark}

\allowdisplaybreaks
\def\be{\begin{equation}}
\def\ee{\end{equation}}

\def\ve{\varepsilon}
\def\vp{\varphi}

\def\arrowk{^\to{\kern -6pt\topsmash k}}
\def\arrowK{^{^\to}{\kern -9pt\topsmash K}}
\def\arrowr{^\to{\kern-6pt\topsmash r}}
\def\arrowvp{^\to{\kern -8pt\topsmash\vp}}
\def\arrowf{^{^\to}{\kern -8pt f}}
\def\arrowg{^{^\to}{\kern -8pt g}}
\def\arrowu{^{^\to}a{\kern-8pt u}}
\def\arrowt{^{^\to}{\kern -6pt t}}
\def\arrowe{^{^\to}{\kern -6pt e}}
\def\tk{\tilde{\kern 1 pt\topsmash k}}
\def\barm{\bar{\kern-.2pt\bar m}}
\def\barN{\bar{\kern-1pt\bar N}}
\def\barA{\, \bar{\kern-3pt \bar A}}

\def\iint{\not \kern-4pt\int}
\begin{document}
\title{A prime number theorem for the majority function}

\author{Jean Bourgain}

\address{School of Mathematics, Institute for Advanced Study, 1
Einstein Drive, Princeton, NJ 08540.}
\email{bourgain\@math.ias.edu}
\thanks{The research was partially supported by NSF grants DMS-0808042 and
DMS-0835373.}
\begin{abstract}
In the paper, the occurrence of zeros and ones in the binary expansion of the primes is studied.
In particular the statement in the title is established.  The proof is unconditional.
\end{abstract}
\maketitle

\section
{Introduction}

Let $N=2^n$ and identify $\{ 0, 1, \ldots, N-1\}$ with $\{0, 1\}^n$ by binary
expansion
$$
x=\sum_{0\leq j< n} x_j 2^j \ \text { with } x_j=0, 1.
$$
Assuming $n$ odd, denote $f:\{0, 1\}^n\to \{0, 1\}$ the majority function.
Motivated by a question of G.~Kalai\cite {Ka}, 
we prove that $f$ does not correlate with the primes, i.e.

\begin{theorem}\label {Theorem1}
Let $\Lambda $ be the Von Mangoldt function. Then
\be
\label{1.1}
\sum_{1\leq x<N} \Lambda (x) f(x) \approx \frac N2.
\ee
\end{theorem}

Note that the majority function is a monotone Boolean function and it was
proven in \cite{B3} that the Moebius function does not correlate with any
monotone Boolean function.
The proof of his property uses the concentration of the Fourier-Walsh
spectrum of monotone Boolean function on `low levels'.
More precisely, expanding
\be
\label{1.2}
f(x) =\sum_{S\subset\{0, \ldots, n-1\}} \hat f(S) w_S(x)
\ee
with
$$
w_S(x) =\prod_{j\in S} \ve_j, \ve_j = 1-2 x_j
$$
the Walsh system on $\{0, 1\}^n$, one exploits that
$$
\sum_{|S|>n^{\frac 12+\ve}} |\hat f(S)|^2
$$
is small for monotone Boolean functions.
This concentration is not  sufficiently strong however to treat $\Lambda$
instead of $\mu$.

Recall that for the majority function, by symmetry,  $\hat f(S) =\hat f(|S|)$ which obey
\be\label{1.3}
|\hat f(k)|^2 \sim \text{\small $\begin{pmatrix} n\\ k\end{pmatrix}^{-1}$} k^{-3/2} \ \text
{ for } \ k>0.
\ee
Hence
\be\label{1.4}
\sum_{|S|=k} |\hat f(S)|^2 \sim k^{-3/2}
\ee
and
\be\label{1.5}
\sum_{|S|>k} |\hat f(S)|^2 \lesssim k^{-1/2}.
\ee
Write
\be\label{1.6}
\sum^N_1 \Lambda(x) f(x) =\frac 12\Big(\sum_1^N\Lambda(x)\Big) +
N\sum_{0<|S|\leq n} \hat\Lambda (S) \hat f(S).
\ee
Introducing some cutoff $n_0<n$, estimate the second term of (2.6) by
\be \label{1.7}
N\sum_{0<|S|\leq n_0} |\hat\Lambda(S)| \ |\hat f(S)|
\ee
+
\be\label{1.8}
N\sum_{n_0<|S|\leq n} |\hat\Lambda(S)| \ |\hat f(S)|.
\ee

Because primes are odd (except for the prime 2), for $S=(1, 0, \ldots, 0)$,
$$
\hat \Lambda (S) =\frac 1N \sum^N_{x=1} \Lambda (x) (1-2x_1) = -\frac
1N\sum^N_{x=1} \Lambda (x) \approx -1.
$$
For $0<|S|<\sqrt n$, $S\not= (1, 0, \ldots, 0)$, it follows from \cite{B2}
that
\be\label{1.9}
|\hat\Lambda(S)|< e^{-c\sqrt n}.
\ee
Taking
\be\label{1.10}
n_0\sim n^{\frac 12-\ve}
\ee
the preceding permits to bound \eqref{1.7} by
\be\label{1.11}
O\Big(\frac N{\sqrt n}+N e^{-c\sqrt n}\sum_{k<n_0} 
\begin{pmatrix} n\\ k\end{pmatrix}^{-\frac 12} k^{-\frac 34}\Big)= O
\Big(\frac N{\sqrt n}\Big).
\ee
On the other hand, if we try to estimate \eqref{1.8} using $L^2$-norm, the
tail  estimate \eqref{1.5} implies
\be\label{1.12}
\eqref{1.8} \leq N\sqrt n\Big(\sum_{|S|> n_0} |\hat f(S)|^2\Big)^{\frac 12}\lesssim N\sqrt n n_0^{-1/4}
\ee
which is not conclusive, no matter how $n_0\ll n$ is chosen.

Hence a more refined analysis is needed, involving more than just the low Fourier-Walsh spectrum of $\Lambda$.

In what follows, we will rely in particular on estimates related to those in the work of Mauduit-Rivat \cite{M-R}, where it was shown that
$\Lambda$ does not correlate with the parity function
\be\label{1.13}
\sigma(x) =e^{i\pi(\sum_{0\leq j< n} x_j)} = w_{\{0, 1, \ldots, n-1\}}(x)
\ee
(rather than the majority function).
See also \cite{B1} from which we will borrow certain estimates.

Before going further, we point out the following easy consequence of \cite{B2} on prescribing binary digits from the primes.

\begin{theorem}\label{Theorem2}
Let $\rho<\frac 47$.
Then, with above notations, taking $r\sim n^\rho$, there are at least $O(2^{-r}\frac Nn)$ primes less than $N$ satisfying
\be\label{1.14}
\sum_1^n x_j >\frac n2+ \frac13 r
\ee
\vfill\eject
and at least $O(2^{-r}\frac Nn)$ primes for which
$$
\sum_1^n x_j< \frac n2-\frac 13 r.
$$
\end{theorem}

It follows indeed from \cite{B2} that for $r<n^{\frac 47-}$, the set
$$
\Omega_1= \{p<N, x_0=x_1=\cdots = x_{r-1}=1\}
$$
satisfies
$$
|\Omega_1| \sim \frac Nn 2^{-r}.
$$
Since also for $1\ll \Delta <\log n$
$$
\Big|\Big\{x<N; x_0 =\cdots= x_{r-1} =1 \text { and } \Big|\sum^{n-1}_{r+1} x_j-\frac {n-r}2\Big| >\Delta \sqrt{n-r}\Big\}
\Big|< e^{-c\Delta^2} N 2^{-r}
$$
necessarily most elements of $\Omega_1$ will satisfy
$$
\Big|\sum^{n-1}_{j=r+1} x_j-\frac {n-r}2\Big| < O(\sqrt {n\log n})
$$
and
$$
\sum^{n-1}_{j=0} x_j >\frac {n+r}2 -O(\sqrt{n\log n}).
$$
The second part of the statement is proven similarly, considering the set
$$
\Omega_0=\{p<N; x_1=\cdots= x_{r-1} =0\}.
$$
Note that it is essential for this argument that $r\gg n^{\frac 12}$.
\vskip.3 true in

\noindent
{\bf Acknowledgement}: The author is grateful to G.~Kalai for bringing up various problems on the digital aspects of arithmetic functions and correspondence on
those results.

\section
{Symmetrization of the Von Mangoldt function}

Returning to the proof of Theorem \ref{Theorem1}, we note that
$$
\sum^N_1\Lambda(x) f(x) \equiv \langle\Lambda, f\rangle =\langle \Lambda_s. f\rangle
$$
where $\Lambda_s$ stands for the symmetrization of $\Lambda$ under the permutation group of $\{0, 1, \ldots, n-1\}$.
Thus
\be\label{2.1}
\Lambda_s =\sum^n_{k=1} \ \frac {\sum_{x\in\Omega_k} \Lambda (x)}{\Big(
{\begin{matrix} n\\ k\end{matrix}}\Big)} 1_{\Omega_k}
\ee
where $\Omega_k=\{ x\in\{0, 1\}^n; \sum x_j=k\}$.

The advantage of introducing $\Lambda_s$ is a reduction of the $L^2$-norm.
For $0\leq \rho\leq 1$, denote $T_\rho$ the usual convolution operator defined by
$$
T_\rho w_S=\rho^{|S|}w_S
$$
and which is a contraction on all $L_p$-spaces.
Write
$$
\langle \Lambda, f\rangle =\langle T_\rho \Lambda, f\rangle +\langle (1-T_\rho)\Lambda_s, f\rangle = (2.2)+ (2.3).
$$
Then
$$
(2.2) =\frac 12 \sum^N_1\Lambda(x) + N\sum_{0<|S|\leq n} \rho^{|S|} \hat\Lambda (S) \hat f(S)\eqno{(2.4)}
$$
and estimate, recalling \eqref{1.9} the second term of (2.4) by
$$
O\Big(\frac N{\sqrt n}\Big) +N\sum_{n_0< k\leq n} \rho^k k^{-3/4}\Big[\sum_{|S|=k} |\hat\Lambda (S)|^2\Big]^{\frac 12}\\
$$
$$
<O\Big(\frac N{\sqrt n}\Big)+ Nn^{\frac 12} \rho^{-n_0}<O\Big(\frac N{\sqrt n}\Big)\qquad\qquad \eqno{(2.5)}
$$
provided, cf \eqref {1.10}, we set
$$
\rho=1-n^{-\frac 12+2\ve}.\eqno{(2.6)}
$$

To estimate (2.3), we decompose further
$$
\Lambda_s=\Lambda_s'+\Lambda^{\prime\prime}_s\eqno{(2.7)}
$$
denoting
$$
\Lambda_s' =\sum_{|k-\frac n2|<\Delta\sqrt n} \ \frac {\sum_{k\in\Omega_k} \Lambda (x)}{\Big(\begin{matrix} n\\ k\end{matrix}\Big)} 
\, 1_{\Omega_k}
$$
with $\Delta\gg 1$ a parameter. Then
$$
|\langle (1-T_\rho)\Lambda_s, f\rangle|\leq |\langle (1-T_\rho)\Lambda_s', f\rangle|+ \Vert\Lambda_s^{\prime\prime} \Vert_1.\eqno{(2.8)}
$$
Estimate
$$
|\langle (1-T_\rho)\Lambda_s', f\rangle|\leq \Vert\Lambda_s'\Vert_2 \ \Vert(1-T_\rho)f\Vert_2
$$
where
$$
\Vert\Lambda_s'\Vert_2= \Bigg\{\sum_{|k-\frac n2|<\Delta\sqrt n} \frac {\big(\sum_{x\in\Omega_k}\Lambda(x)\big)^2}{\Big({\begin{matrix} n\\ k
\end{matrix}}\Big)} \Bigg\}^{\frac 12}
\leq
$$
$$
\sqrt N\Bigg[\max_{|k-\frac n2|<\Delta\sqrt n} \ \frac {\sum_{x\in\Omega_k} \Lambda (x)}{\Big({\begin{matrix} n\\ k\end{matrix}}\Big)}\Bigg]^{\frac 12}\lesssim
$$
$$
n^{\frac 14} e^{C\Delta^2} \Big[\max_k\sum_{x\in\Omega_k} \Lambda(x)\Big]^{\frac 12}\eqno{(2.9)}
$$
and, again form \eqref{1.3}, (2.6)
$$
\begin{aligned}
\Vert (1- T_\rho)f\Vert_2 &\leq \sqrt N\Big[\sum_k (1-\rho^k)^2 k^{-3/2}\Big]^{\frac 12}\\
&\leq \sqrt N\Big[ n_0^{-1/2} +\sum_{k\leq n_0} k^{1/2} (1-\rho)^2 \Big]^{\frac 12} \lesssim n^{-\frac 18+2\ve} \sqrt N.
\end{aligned}
\eqno{(2.10)}
$$
Hence
$$
|\langle (1-T_\rho)\Lambda_s', f\rangle |\lesssim n^{\frac 18 +2\ve} e^{c\Delta^2}\Big\{\max_k\Big[\frac 1N\sum_{x\in\Omega_k} \Lambda (x) \Big]\Big\}^{\frac 12}
N.\eqno{(2.11)}
$$
Next
$$
\Vert \Lambda_s^{\prime\prime} \Vert_1 =\sum_{|k-\frac n2|\geq \Delta\sqrt n}\Lambda (x).
$$
Let $R\in \mathbb Z_+, R<\log n$ and estimate, again using the correlation estimates of $\Lambda$ with low order Walsh functions
$$
\begin{aligned}
&\sum^N_1\Lambda(x) \Big|\frac n2 -\sum x_j\Big|^{2R} \leq \sum^N_1 \Lambda (x) \Big|\sum_0^{n-1} \ve_j\Big|^{2R}\\
&\lesssim (CR)^R n^RN+ (CR)^R\Big(\sum_{o<|S|\leq 2R} |\hat\Lambda (S)|\Big) \lesssim (CR)^R n^RN.
\end{aligned}\eqno{(2.12)}
$$
Therefore
$$
\sum_{|\frac n2 -\sum x_j|>\Delta\sqrt n} \Lambda(x) < e^{-c\Delta^2}N.\eqno{(2.13)}
$$
It remains to establish a bound on
$$
\sum_{x\in\Omega_k} \Lambda(x)\eqno{(2.14)}
$$
for $|k-\frac n2|\leq \Delta \sqrt n$ in (2.11).

\section
{Distribution of the sum of the digits of the primes}

Our remaining task is to bound (2.14) in the range $k=\frac n2 +O(\sqrt n)$.
Take a bumpfunction $\eta$ on $\mathbb R$ s.t. $\hat\eta\geq 0, \hat\eta (0)=1$ and $\text{supp\,} \eta\subset [-\frac 12, \frac 12]$ say.

Clearly
\be\label{3.1}
\sum_{x\in\Omega_k} \Lambda(x) =\int_{-\frac 12}^{\frac 12} \Big[\sum^N_1 \Lambda (x) e^{i\lambda(\sum^n_1 x_j -k)}\Big]\eta(\lambda) d\lambda
\ee
and we evaluate
\be\label{3.2}
\sum^N_1\Lambda (x) U_\lambda(x)
\ee
where
\be\label{3.3}
U_\lambda(x) =e^{i\lambda(\sum\limits_0^{n-1} x_j)}.
\ee
This issue is very similar to the case of the Morse function $(\lambda=\pi)$ considered by Mauduit-Rivat in \cite{M-R}.
Thus we will use the Vinogradov type I-II sum approach from \cite{M-R}.
In what follows, we will in fact rely on the presentation in \cite{B1} (where the Mo\"ebius function rather than $\Lambda$ is considered, but there is no
essential difference here between these cases.)

The Fourier coefficients of $U_\lambda$ obey an estimate
\be\label{3.4}
\Big|\hat U_\lambda(k)\Big|\lesssim e^{-c\lambda^2 n}.
\ee 
The argument is similar to Lemma 2 in \cite{B1}.
In case of the Morse  sequence $w_{\{0, 1, \ldots, n-1\}} = U_\pi$, one has in particular $\Vert\hat U_\pi\Vert_\infty < e^{-cn}$ which is stronger than
\eqref{3.4} for small $\lambda$. This is the most significant difference compared with \cite{B1}.

Recall some terminology.
Let $n=m_1+m_2, M_1=2^{m_1}, M_2=2^{m_2}, m_1\leq m_2$.

Type-$II$ sums are of the form
\be\label{3.5}
\sum_{\substack {x^1\sim M_1\\ x^2\sim M_2}}  a_{x_1} b_{x_2} U_\lambda (x^1.x^2)
\ee
where $a_{x_1}, b_{x_2}$ are (arbitrary) bounded sequences (in fact obtained) as multiplicative convolutions of $\Lambda$ and $\mu$) and we may assume 
$M_1>N^{\frac 13}$.
For the Type-$I$ sums, we set $b_{x_2} =1$.
Of course, the analysis of Type-$II$ sums applies equally well to the Type-I sum but for the latter, also other considerations will be involved
when $M_1$ is small.

We start by recalling the Type-$II$ bound (2.31) from \cite{B1}, which in view of \eqref{3.4} becomes
\be\label{3.6}
|\eqref{3.5}|\lesssim N(L^{-c_1} +L^2 M_1^{-c_2}+L^{C _3} M_1^{-c\lambda^2})
\ee 
where $c_1, c_2, C_3$ are some constants, $L$ a parameter (note that \cite{B1} treats the case of an arbitrary Walsh function $w_S$, while for our purpose only
the case $S=\{0, 1, \ldots, n-1\}$ is of relevance).

Optimizing \eqref{3.6} in $L$ gives a bound of the form
\be\label{3.7}
NM_1^{-c'\lambda^2}.
\ee

Next, according to \cite{B1}, (3.2') and \eqref{3.4}, the following estimate on Type-I sums is gotten
\be\label{3.8}
M_1^2M_2\Vert\hat U_\lambda\Vert_\infty \lesssim NM_1 \, e^{-c\lambda^2n}.
\ee
Assuming $M_1>N^{\frac 13}$, \eqref{3.7} gives a bound $Ne^{-c\lambda^2n}$ on Type-$II$ sums.
The Type-$I$ sums may be estimated using either \eqref{3.7} or \eqref{3.8}, hence satisfy a bound $N. e^{-c\lambda^4n}$, which is conclusive provided
\be\label{3.9}
\lambda> n^{-\frac 14+\ve}.
\ee
The range \eqref{3.9} is not quite sufficient for our needs.
Consequently assume
\be\label{3.10}
n^{-\frac 12+\ve}< \lambda< n^{-\frac 14+\ve}
\ee
and in view of the already available estimates \eqref{3.7}, \eqref {3.8}, also
\be\label{3.11}
c\lambda^2 n\lesssim m_1< n^\ve \lambda^{-2}.
\ee
Take
\be
\label{3.12}
m_1\ll m\ll n
\ee
to specify and decompose
$$
x=(y, z)\in\{0, 1\}^m\times \{0, 1\}^{n-m}.
$$
Write
$$
U_\lambda (x) = e^{i\lambda(\sum^{m-1}_{0} y_j)} \ e^{i\lambda(\sum^{n-1}_{m} z_j)} = U(y)V(z).
$$
Hence
\be
\label{3.13}
\sum_{x^1\sim M_1} \Big|\sum_{x^2\sim M_2} U_\lambda (x^1.x^2)\Big| \leq \sum_z\ \sum_{x^1\sim M_1} \Big|\sum_{y\equiv -z(\text{mod\,} x^1)} U(y)\Big|.
\ee

Some further manipulation of $U(y)$ is needed.
Write $\ve_j= 1-2y_j$ and
\be
\label{3.14}
U(y) = e^{i\frac\lambda 2m} \Big(\cos \frac \lambda 2\Big)^m \ \prod^m_{j=1} (1+i\ve_j tg \, \frac \lambda 2).
\ee
Expanding the last factor of \eqref{3.14} in the Walsh system
$$
\begin{aligned}
\prod^m_{j=1} \Big(1+ i\ve_j tg \, \frac \lambda 2\Big) &=\sum_{k\leq k_1} (itg\,\frac\lambda 2)^k \sum_{|S|=k} w_S(\ve)+\sum_{k_1< k\leq m}\cdots\\
&=(3.15)+(3.16).
\end{aligned}
$$
Taking
$$
k_1\sim \lambda^2 m\eqno{(3.17)}
$$
gives
$$
\Vert(3.16)\Vert_2< \sum_{m\geq k>k_1}|\lambda|^k \Big(\begin{matrix} m\\ k\end{matrix}\
\Big)^{\frac 12} < 1
$$
and the contribution of (3.16) in \eqref{3.13} is bounded by
$$
\Big( \cos\frac \lambda 2\Big)^m 2^m 2^{n-m} < e^{-c\lambda^2m}N.\eqno{(3.18)}
$$
Next, in (3.15), expand
$$
h=\sum_{|S|=k} w_S
$$
in a regular Fourier series
$$
h(y) =\sum_{r=0}^{2^m-1} \hat h(r) e\Big(\frac {ry}{2^m}\Big).\eqno{(3.19)}
$$
Fixing $0\leq r< 2^m$, substituting in \eqref {3.13}, we obtain
$$
\begin{aligned}
&\sum_{x^1\sim M_1} \Big|\sum_{y\equiv -z (\text{mod\,} x^1)} \ e\Big(\frac {ry}{2^m}\Big)\Big|\lesssim\\
&\frac {2^m}{M_1} \sum_{x^1\sim M_1} \ 1_{[\Vert\frac {x^1r}{2^m}\Vert<n\frac {M_1}{2^m}]}.
\end{aligned}\eqno{(3.20)}
$$
Let $\delta>0$ be another parameter and assume that
$$
\sum_{x^1\sim M_1} 1_{[\Vert\frac {x^1r}{2^m}\Vert <n \frac {M_1}{2^m}]}>\delta M_1.\eqno{(3.21)}
$$

By the pigeonhole principle, there is some $q'\lesssim \frac 1\delta$ s.t. $\Vert\frac{q'r}{2^m}\Vert < n\frac {M_1}{2^m}$ and therefore we get
$$
\frac r{2^m} =\frac aq+\theta\eqno{(3.22)} 
$$
with
$$
q\lesssim \frac 1\delta, (a, q)=1 \ \text { and } \ |\theta|<\frac {nM_1}{q{2^m}}.\eqno{(3.23)}
$$
Assuming
$$
\delta> 2^{-\frac m2}\eqno{(3.24)}
$$
it follows that $\Vert\frac {x^1 r}{2^m}\Vert \geq \Vert\frac {x^1a}q\Vert-\frac 1{2^{m-2m_1}} >\frac\delta 2$, unless
$x^1 a\equiv 0 (\text{mod\,} q)$.
If $x^1a\equiv 0 (\text{mod\,} q), \Vert\frac {x^1r}{2^m}\Vert = x^1|\theta|$
and we obtain the condition
$$
x^1<\frac {nM_1}{2^m|\theta|}.
$$
In view of (3.25), this implies that
$$
|\theta|\lesssim \frac n{2^m\delta}.\eqno{(3.25)}
$$
Let $\frac r{2^m}$ satisfy (3.22) with
$$
q< \frac 1\delta \ \text { and } \ |\theta|\lesssim \frac n{2^m\delta}.
\eqno{(3.25)}
$$
We estimate $\hat w_S(r)$.
Thus, letting $\vp =\frac r{2^m}$
$$
\begin{aligned}
\hat w_S(r) &= 2^{-m} \sum_{(x_0, \ldots, x_{m-1})\in \{0, 1\}^m}
e^{2\pi i\vp (\sum_{j=0}^{m-1} 2^j x_j)+i\pi \sum_{j\in S} x_j}\\
&= 2^{-m} \prod_{j\not\in S}(1+ e^{2\pi i 2^j\vp}) \prod_{j\in S} (1- e^{2\pi i 2^j\vp})
\end{aligned}
$$
and
$$
\begin{aligned}
|\hat w_S(r)|& =\prod_{j\not\in S} |\cos \pi 2^j\vp|\prod_{j\in S} |\sin \pi 2^j\vp|\\
&\leq \prod_{\substack{j\not\in S\\ j< m-J}} \Big(\Big|\cos 2\pi 2^j\frac aq\Big|+ \pi 2^j|\theta|\Big) \ \prod_{\substack{j\in S\\ j< m-J}}
\Big(\Big|\sin \pi 2^j\frac aq\Big|+\pi 2^j|\theta|\Big)\end{aligned}
\eqno{(3.27)}
$$
with $1\ll J\ll m$ to specify.

By (3.26), $2^j|\theta| \lesssim \frac n\delta 2^{-(m-j)}\lesssim \frac n\delta. 2^{-J}$ and we take
$$
J\sim \log\frac 1\delta+k\eqno{(3.28)}
$$
as to ensure that
$$
|\hat w_S(r)| \leq\prod_{\substack {j\notin S\\ j< m-J}} \Big|\cos 2\pi 2^j\frac aq\Big| \ \prod_{\substack{j\in S\\ j< m-J}} 
\Big|\sin \pi 2^j\frac aq\Big|+\delta.
\eqno{(3.29)}
$$
Recall that $|S|=k\leq k_1\sim \lambda^2m$.
It follows that there is an interval $\{j_0, \ldots, j_1-1\}$ in $\{0, \ldots, [\frac m2]\}$ of size
$$
j_1-j_0>\frac m{2k_1}\eqno{(3.30)}
$$
which is disjoint from $S$.
The first factor in (3.29) is then majorized by
$$
\begin{aligned} 
\prod_{j\in I} 
\Big|\cos 2\pi 2^j\frac aq\Big| &=\frac 1{2^{j_1-j_0}} \Big|\sum_{u=0}^{2^{j_1-j_0} -1} e^{2\pi i  2^{j_0}\frac aq}\Big|\\
&\leq\frac q{2^{j_1-j_0}}< \frac 1{\delta 2^{\frac m{2k_1}}}
\end{aligned}
\eqno{(3.31)}
$$
provided $q$ is not a power of 2.
On the other hand, if $q$ is a power of 2, then $\sin \pi 2^j\frac aq = 0$ for $j\gtrsim \log \frac 1\delta$ and we conclude that
$$
|\hat w_S(r)|< \delta+\frac 1{\delta 2^{\frac m{2{k_1}}}} < \delta+\frac 1\delta e^{-c\lambda^{-2}}\eqno{(3.32)}
$$
except if $S\subset \{0, 1, \ldots, J\} \cup \{m-J, \ldots, m-1\}$.

Consequently, the contribution of the $k$-term of (3.15) in \eqref{3.13} may be estimated as follows
$$
2^n\Big(\cos \frac \lambda 2\Big)^m \Big|tg \frac \lambda 2\Big|^k 
\Big\{\Vert\hat h\Vert_1\, \delta+\begin{pmatrix} m \\ k\end{pmatrix} \Big(\delta+\frac 1\delta\,
e^{-c\lambda^{-2}}\Big) +\begin{pmatrix} 2J\\ k\end{pmatrix} \max_{|S|=k} \Vert\hat w_S\Vert_1\Big\}.
\eqno{(3.33)}
$$
with $J$ given by (3.28).

Making a suitable approximation of the step-function by Fourier-truncation (cf. \cite{B1} for details), with an $L^1$-error at
most $m^{-k}$ say, we ensure that
$$
\Vert\hat w_S\Vert_1 < (ck \log n)^k\eqno{(3.34)}
$$
and hence
$$
\Vert\hat h\Vert_1< \begin{pmatrix} m\\ k\end{pmatrix} (ck\log n)^k.\eqno{(3.35)}
$$ 
Substituting (3.34), (3.35) in (3.33), we find
$$
(3.33)< 2^n e^{-c\lambda^2m}
\Big\{ m^{2k}\delta+m^k \delta^{-1} e^{-c\lambda^{-2}}+ \Big( 1+ \frac {c\log \frac 1\delta} k\Big)^k
(ck\lambda\log n)^k\Big\}.\eqno{(3.36)}
$$
Taking $\delta =m^{-2k}$ gives
$$
\begin{aligned}
(3.33)&< 2^n \, e^{-c\lambda^2m}\Big(1+m^{3k_1} e^{-c\lambda^{-2}}+(Ck_1(\log n)^2 \lambda )^k\Big)
\qquad \qquad\\
&< 2^n \, e^{-c\lambda^2m} \Big(1+ e^{c(\log n)\lambda^2m-c\lambda^{-2}}+ (C(\log n)^2 \lambda^3 m)^k\Big).
\qquad
\end{aligned}
\eqno{(3.37)}
$$
Recalling \eqref{3.10}-\eqref{3.12}, take
$$
m=\frac c{(\log n)^2} \min (\lambda^{-3}, n).\eqno{(3.38)}
$$
Then
$$
(3.37) < 2^n e^{-c\lambda^2m}< 2^n e^{-c(\log n)^{-2} \min (\lambda^{-1}, \lambda^2n)}\eqno{(3.39)}
$$
which gives a bound for the (3.15)-contribution to \eqref{3.13}.

Thus we proved that if $n^{-\frac 12+\ve}<\lambda< n^{-\frac 14+\ve}$ and $\lambda^2 n\lesssim m_1<\frac {n^\ve}{\lambda^2}$, then
$$
\eqref{3.13}< N e^{-n^\ve}.\eqno{(3.40)}
$$
Summarizing, we conclude that \eqref{3.2} may certainly be bounded by $\frac Nn$ provided $\lambda> n^{-\frac 12+\ve}$.

Consequently, substituting in \eqref{3.1} gives
$$
\sum_{x\in\Omega_k} \Lambda(x) < N n^{-\frac 12+\ve}.\eqno{(3.41)}
$$

\section
{Conclusion of the proof of Theorem 1}

Substitution of (3.41) in (2.9) gives
\be\label{4.1}
\Vert\Lambda_s'\Vert_2< e^{C\Delta^2}n^\ve N
\ee
and in (2.11)
\be\label{4.2}
|\langle (1-T_\rho)\Lambda_s', f\rangle| < e^{C\Delta^2} n^{-\frac 18+3\ve} N.
\ee
Recalling (2.5) and (2.13), we proved that
\be\label{4.3}
\begin{aligned}
\langle \Lambda, f\rangle&< O\Big( \frac N{\sqrt n}\Big) +e^{C\Delta^2} n^{-\frac 18 +3\ve} N+ e^{-c\Delta ^2}N\\
&< n^{-c} N
\end{aligned}
\ee
for some constant $c>0$, by suitable choice of $\Delta$.

Hence, Theorem \ref{Theorem1} holds in the more precise form
\be
\sum^N_1\Lambda(x) f(x) =\frac N2+O(n^{-c} N).
\ee


\begin{thebibliography}
{xxxxxxxx}

\bibitem[B1]{B1} J.~Bourgain,
\emph{Moebius-Walsh correlation bounds and an estimate from Mauduit and Rivat}, to appear in J.~Analyse.

\bibitem[B2]{B2} J.~Bourgain, \emph{Prescribing the binary digits of the primes}, to appear in Israel J.~Math.

\bibitem[B3]{B3} J.~Bourgain, \emph{On the Fourier-Walsh spectrum of the Moebius function}, to appear in Israel J. Math.

\bibitem[Ka] {Ka} G.~Kalai, \emph{Private communications}

\bibitem [M-R] {M-R} C.~Mauduit, J.~Rivat, \emph{Sur un probl\`eme de Gelfand; la sommes des chiffres des nombres premiers},
Annals of Math {\bf 171} (2010), 1591--1646.







\end{thebibliography}
\end{document}